\documentclass[final,3p]{elsarticle}
\usepackage{amsmath,amssymb,amsthm}

\newcommand{\drm}{\mathrm{d}}

\newcommand{\PP}{\vecb{P}}
\newcommand{\HH}{\vecb{H}}

\newcommand{\LL}{\vecb{L}}

\newcommand{\VV}{\vecb{V}}
\newcommand{\Vz}{\vecb{V}^0}

\newcommand{\helm}{\mathbb{P}}
\newcommand{\stokes}{\mathbb{S}}

\newcommand{\vecb}[1]{\mathbf{#1}}
\newcommand{\ee}{\vecb{e}}
\newcommand{\ff}{\vecb{f}}
\newcommand{\uu}{\vecb{u}}
\newcommand{\vv}{\vecb{v}}
\newcommand{\ww}{\vecb{w}}
\newcommand{\xx}{\vecb{x}}

\newtheorem{thm}{Theorem}
\newtheorem{rem}{Remark}

\journal{Journal of Computational Physics}

\begin{document}

\begin{frontmatter}

\title{Pressure-induced locking
in mixed methods for time-dependent (Navier--)Stokes equations}

\author{Alexander Linke\corref{cor1}%
   }
\ead{ alexander.linke@wias-berlin.de}
 \address{Weierstrass Institute, Mohrenstr. 39, 10117 Berlin, Germany}

\author{Leo G. Rebholz 
}
\ead{ rebholz@clemson.edu }
 \address{Department of Mathematical Sciences, Clemson University, Clemson SC 29634, USA}

%
%


\end{frontmatter}

\section{Introduction}
We consider inf-sup stable mixed methods for the time-dependent incompressible Stokes and Navier--Stokes equations,
extending earlier work on the steady (Navier--)Stokes problem
\cite{JLMNR17}. A locking phenomenon is identified for classical inf-sup stable methods like the Taylor--Hood or the Crouzeix--Raviart elements by a novel, elegant and simple numerical analysis and corresponding numerical experiments, whenever the momentum balance is dominated by forces of a gradient type.
More precisely, a reduction of the $\LL^2$ convergence order for high order methods, and even a complete stall of the $\LL^2$ convergence order for lowest-order methods on preasymptotic meshes is predicted by the analysis and practically observed. On the other hand, it is also shown that (structure-preserving) {\em pressure-robust} mixed methods do not suffer from this locking phenomenon, even if they are of lowest-order.

The short note contributes to the recent scholarly debate on the accuracy of low-order structure-preserving space discretizations, e.g. with respect to the treatment of gradient fields in the momentum balance by well-balanced schemes for the shallow water or compressible Euler equations, and the accuracy of (non-structure-preserving) high-order space discretizations \cite{HAT18}.
It demonstrates that the structure-preserving,
well-balanced property can be achieved in our setting, if certain
discretely divergence-free velocity test functions are even
weakly divergence-free in the sense of $\LL^2$ \cite{JLMNR17}
--- without needing to know the exact form of the equilibrium solution, which is a typical disadvantage of well-balanced
schemes for hyperbolic conservation laws \cite{HAT18}.
Thus, the short note builds a bridge between inf-sup stable mixed finite elements for the Navier--Stokes (NS) equations and well-balanced schemes for hyperbolic conservation laws, which have traditionally not too much exchange of knowledge.

\section{Time dependent Stokes and the finite element space discretization}
Consider the time-dependent incompressible Stokes equations
with homogoneous Dirichlet boundary conditions in $(0, T] \times \Omega$: find $(\uu,p)$ satisfying
\begin{equation} \label{eq:cont:problem}
\begin{split}
  \uu_t - \nu \Delta \uu + \nabla p & = \ff, \\
  \nabla \cdot \uu & = 0, 
\end{split}
\end{equation}
with pressure assumed to be normalized, and the divergence-free initial value is prescribed as $\uu(0, \xx) = \uu_0(\xx)$.
For simplicity, the domain $\Omega$ is assumed
as convex to guarantee elliptic regularity. Further, the discussion
is restricted to smooth solutions $(\uu, p)$.

We denote the $L^2(\Omega)$
inner product by $(\cdot,\cdot)$, the $H^k(\Omega)$ norm by  $\| \cdot \|_k$, and assume that conforming 
finite element velocity-pressure spaces $(\VV_h,Q_h)$ satisfy the Babuska--Brezzi condition \cite{brezzi:fortin} (extension of our analysis to stable nonconforming methods is straight-forward but requires significant extra notation).
The discretely divergence free velocity space is defined by $\VV_h^0:=\{\vv_h \in \VV_h,\ (\nabla \cdot \vv_h,q_h)=0 \ \forall q_h\in Q_h \}.$

The standard finite element spatial discretization
on shape-regular triangulations
is given as follows: for all $t \in (0, T]$ search for $(\uu_h, p_h) \in (\VV_h, Q_h)$ such that 
\begin{equation} \label{eq:disc:problem}
\begin{split}
  (\dot{\uu}_h, \vv_h) - \nu (\nabla \uu_h, \nabla \vv_h)
    -(p_h, \nabla \cdot \vv_h)  + (\nabla \cdot \uu_h, q_h) = (\ff, \vv_h), \\
\end{split}
\end{equation}
for all $(\vv_h, q_h) \in (\VV_h, Q_h)$.  The discrete initial value is prescribed
as $\uu_h(0) := \helm_h(\uu_0)$, where $\helm_h$ is a
discrete Helmholtz--Hodge projection \cite{JLMNR17} into the discretely divergence-free space, defined by: Given $\ww \in L^2(\Omega)$, $\helm_h(\ww)\in \VV_h^0$ satisfies
\[
(\helm_h(\ww),\vv_h)=(\ww,\vv_h) \ \forall \vv_h \in \VV_h^0.
\]
We will also utilize a $\HH^1_0(\Omega)$ projection onto $\VV_h^0$, which is called the discrete Stokes projection, and is denoted $\stokes_h$ and defined by: Given $\ww \in \HH^1(\Omega)$, find $\stokes_h(\ww)\in \VV_h^0$ satisfying
\[
(\nabla \stokes_h(\ww),\nabla \vv_h) = (\nabla \ww,\nabla \vv_h) \ \forall \vv_h \in \VV_h^0.
\]
We note that due to the elliptic regularity, i.e., the convexity
of the domain $\Omega$,
and due to the Babuska--Brezzi condition
both $\helm_h$ and $\stokes_h$
have optimal approximation properties 
on divergence-free vector fields 
in both the $\LL^2$ and the $\HH^1$ norms \cite{ALM18}.

\section{A new a-priori error analysis for flows with gradient-dominated momentum balances}

We now present a new a-priori error analysis that reveals precisely how locking and suboptimal convergence 
can occur in flows where gradient forces dominate the momentum balance (e.g. when $\nu \ll 1$).  To begin the analysis, for the discrete velocity solution $\uu_h$, we make the ansatz
$\uu_h := \ee_h + \helm_h(\uu)$, since we will derive
a supercloseness result. Note that it holds $\ee_h \in \Vz_h$.
Testing \eqref{eq:disc:problem} by $\ee_h$ yields
\begin{equation*}
\begin{split}
  (\dot{\uu}_h, \ee_h) + \nu (\nabla \uu_h, \nabla \ee_h)
    & = (\ff, \ee_h) \\
    & = (\uu_t - \nu \Delta \uu + \nabla p, \ee_h) \\
    & = (\helm_h(\uu_t), \ee_h) + \nu (\nabla \stokes_h(\uu), \nabla \ee_h)
      + (\helm_h(\nabla p), \ee_h).
\end{split}
\end{equation*}
Exploiting
$
  (\helm_h(\uu_t), \ee_h) =
   \left (\frac{\drm}{\mathrm{dt}} \helm_h(\uu), \ee_h \right ),
$
and $\uu_h = \ee_h + \helm_h(\uu)$, we obtain
\begin{equation*}
\begin{split}
  \frac{1}{2} \frac{\drm}{\mathrm{dt}} \| \ee_h \|_0^2
  + \nu \| \nabla \ee_h \|_0^2
   &  = 
     \nu (\nabla (\stokes_h(\uu) - \helm_h(\uu)), \nabla \ee_h)
     + (\helm_h(\nabla p), \ee_h).
\end{split}
\end{equation*}
Now using Cauchy--Schwarz and Young inequalities for
the first right hand side term, we estimate
\begin{equation}
  \frac{\drm}{\mathrm{dt}} \| \ee_h \|_0^2 + \nu \| \nabla \ee_h \|_0^2
   \leq \nu \| \nabla (\stokes_h(\uu) - \helm_h(\uu)) \|_0^2
    + 2 (\helm_h(\nabla p), \ee_h). \label{error1}
\end{equation}
It is fundamental to observe that pressure-robust and classical mixed methods differ
qualitatively in how the term $(\nabla p, \ee_h)$ can be treated.
Even though
it holds for the continuous Helmholtz--Hodge projector
$\helm(\nabla p)=0$, i.e.,
the divergence-free part of $\nabla p$ vanishes exactly
\cite{JLMNR17}, the expression $(\nabla p, \ee_h)$
may represent a certain consistency error
of an appropriate discrete Helmholtz--Hodge projector
for non-pressure-robust (i.e., non structure-preserving) space discretizations \cite{JLMNR17}.
Since $\nabla p$ balances the sum of all gradient parts
in $\ff - \uu_t + \nu \Delta \uu$ in the sense of the
Helmholtz--Hodge
decomposition, different behaviors of different space
discretizations reflect their ability to deal with
dominant gradient fields in the momentum balance, bulding a connection to certain well-balanced schemes
for (vector-valued) hyperbolic conservation laws
\cite{HAT18}. Note that
$\| \nabla (\stokes_h(\uu) - \helm_h(\uu)) \|_{\LL^2}$ converges to $0$ with the optimal rate for an $\HH^1$ norm.

We consider below the two cases separately: the pressure-robust case (here, divergence-free Scott--Vogelius elements) 
for which it holds
\begin{equation}
  (\nabla p, \ee_h) = -(p, \nabla \cdot \ee_h) = 0, \label{robust1}
\end{equation}
and the non-pressure-robust case.
 We consider the pressure-robust case first, and combining \eqref{robust1} with \eqref{error1} immediately implies the following result.

\begin{thm} \label{thm:pr}
For conforming, pressure-robust, inf-sup stable
space discretizations \eqref{eq:disc:problem} of
\eqref{eq:cont:problem}, it holds for all $T > 0$
$$
  \| \ee_h(T) \|_0^2 + \nu \| \nabla \ee_h \|^2_{L^2((0, T); \LL_2)}
    \leq \nu \| \nabla (\stokes(\uu) - \helm_h(\uu)) \|^2_{L^2((0, T); \LL_2)}.
$$
\end{thm}
\begin{rem}
Theorem \ref{thm:pr} reveals a remarkable robustness
of pressure-robust space discretizations with respect to
small viscosities $\nu \ll 1$. Indeed, for $\nu \to 0$
and for all $0 < t < T$
one obtains that $\uu_h(t) \to \helm_h(\uu)(t)$, i.e.,
for smaller and smaller viscosities,
$\uu_h(t)$ converges to the (discretely divergence-free) best approximation of $\uu(t, \cdot)$
in the $\LL^2$ sense, yielding optimal $\LL^2$ convergence
on preasymptotic meshes. 
Moreover, for fixed $\nu$  the error
$\| \nabla \ee_h \|^2_{L^2((0, T); \LL_2)}$ will converge
optimally on resolved meshes,
leading to optimal $\LL^2$ convergence by duality.
\end{rem}

For the case of non-pressure-robust inf-sup stable discretizations, which includes the Taylor--Hood element,
the term $(\helm_h(\nabla p), \ee_h)$ can only be estimated. Standard estimates for the time-dependent
Stokes problem
apply a discrete $\HH^{-1}$ estimate
$$
  (\nabla p, \ee_h) = -(p, \nabla \cdot \ee_h)
    = -(p - \pi_h(p), \nabla \cdot \ee_h)
    \leq \| p - \pi_h(p) \|_0 \cdot \| \nabla \ee_h \|_0,
$$
where $\pi_h(p)$ denotes the
$L^2$ best approximation of the pressure $p$ in the discrete pressure space. While this term goes to zero with the optimal
(pressure) convergence rate in $L^2$, one can reasonably bound 
this term in the time-dependent setting only by something like
\begin{equation} \label{eq:helmh:standard:consist:err}
  \| p - \pi_h(p) \|_0 \cdot \| \nabla \ee_h \|_0
    \leq \frac{1}{\nu} \| p - \pi_h(p) \|_0^2 + \nu \| \nabla \ee_h \|_0^2,
\end{equation}
in order to hide it in the left hand side of \eqref{error1}. However, such a
standard estimate can be terribly pessimistic for small viscosities
$\nu \ll 1$.

In order to derive a sharper estimate in the case of finite
time intervals $(0, T]$ and small viscosities $\nu$, we will now
estimate the term $(\helm_h(\nabla p), \ee_h)$ directly in
$\LL^2$. Therefore, we will assume that the discrete pressure $Q_h$
space contains a $H^1$-conforming subspace, which is elementwise
at least affine. Denoting the Langrange interpolation in this discrete pressure space by $L_h$, it holds
$$
  (\nabla p, \ee_h) = (\nabla p - \nabla L_h p, \ee_h)
    \leq \| \nabla (p - L_h p) \|_0 \cdot \| \ee_h \|_0,
$$
since $\ee_h$ is discretely divergence-free. Now, this
term can be estimated by
$$
  (\nabla p, \ee_h) \leq \frac{T}{2} \| \nabla p - \nabla (L_h p) \|_0^2
    + \frac{1}{2 T} \| \ee_h \|_0^2.
$$
Combining this with the above estimates, we obtain
\begin{equation}
  \frac{\drm}{\mathrm{dt}} \| \ee_h \|_0^2 + \nu \| \nabla \ee_h \|_0^2
   \leq \nu \| \nabla (\stokes_h(\uu) - \helm_h(\uu)) \|_0^2
    + T \| \nabla (p - L_h p) \|_0^2
    + \frac{1}{T} \| \ee_h \|_0^2,
\end{equation}
which is amenable for the Gronwall inequality, because only the (harmless) exponential term 
$
  \exp^{\int_0^t \frac{1}{T} \, \mathrm{ds}} \leq e$ 
  will arise from an application of the inequality.  We have proven the following result.
\begin{thm} \label{thm:cl}
For conforming, inf-sup stable
space discretizations \eqref{eq:disc:problem} of
\eqref{eq:cont:problem}, it holds for all $T > 0$
$$
  \| \ee_h(T) \|_0^2 + \nu \| \nabla \ee_h \|^2_{L^2((0, T); \LL_2)}
    \leq e \nu \| \nabla (\stokes(\uu) - \helm_h(\uu)) \|^2_{L^2((0, T); \LL_2)} + e T \| \nabla (p - L_h p) \|_0^2.
$$
\end{thm}
\begin{rem}
For small $\nu \ll 1$ and
a fixed time interval $(0, T]$ one gets now
an $\LL^2$ convergence order for the discrete velocities
equal to the approximation order of the discrete pressure space (or appropriate subspace)
in the $H^1$ (!) norm.
i) Therefore, one does not get any convergence
order for elements with $P_0$ discrete pressures such as in 
the Bernardi--Raugel or Crouzeix--Raviart finite element methods. Then,
the classical estimate \eqref{eq:helmh:standard:consist:err}
shows merely some asymptotic convergence rates for very fine meshes.
ii) For the Taylor--Hood element
this estimate predicts a (suboptimal)
first-order convergence in the $\LL^2$ norm, losing two orders of
convergence.
iii) For the
mini element one loses one order of convergence in $\LL^2$,
since it approximates velocities with first
order in the $\HH^1$ norm and the discrete pressures with
second order in the $L^2$ norm.
\end{rem}
Thus, classical (i.e. non-pressure-robust) inf-sup stable mixed methods
for incompressible flows require high-order discrete
pressure (!) approximations, in order to get
accurate (although still suboptimal) discrete velocities,
since the discrete
Helmholtz--Hodge projector $\helm_h(\nabla p)$ of classical mixed methods does not exactly vanish
and couples pressure and velocity errors via the pressure-dependent (!) definition of the space of discretely divergence-free vector field $\VV^0_h$.
Similarly, the authors of \cite{HAT18}
argue that
well-balanced schemes allow to reduce the approximation order
of the space discretization in hyperbolic conservation laws.


\section{Numerical Experiments}

We give results here for two numerical tests: time dependent Stokes approximation of a problem with known analytical solution, and time dependent NS approximation of the Chorin vortex decay problem.  In both tests we use small and large viscosity, and varying element choices.  For $\nu=1$ we observe the expected optimal convergence, but when $\nu\ll 1$ we observe precisely the behavior predicted by our (time dependent Stokes) analysis: pressure-robust methods converge optimally, while non-pressure-robust methods converge suboptimally or even lock.

\subsection{Suboptimal convergence and locking when $\nu\ll 1$}

\begin{table}[h!]
{\footnotesize
			\centering
			\begin{tabular}{|c|c | c | c || c | c || c | c || c | c |}
				\hline
				& & \multicolumn{2}{c ||}{TH $(P_2,P_1)$} & \multicolumn{2}{c || }{SV $(P_2,P_1^{disc})$} 
				& \multicolumn{2}{c ||}{Mini $(P_1^b,P_1)$} & \multicolumn{2}{c | }{CR $(P_1^{nc},P_0)$} 
				\\ \hline 
				$\nu$ & $h$ & $\| (u- u_h)(T) \|_0 $ & Rate & $\| (u- u_h)(T) \|_0 $ & Rate & $\| (u- u_h)(T) \|_0 $ & Rate & $\| (u- u_h)(T) \|_0 $ & Rate  \\ \hline
				1 & 1/8 & 1.260e-4 & -      & 9.064e-5 & - & 5.655e-3    & -    &      3.039e-3 & - \\ \hline
				1 & 1/16 & 1.532e-5 & 3.00 & 1.134e-5 &  3.00 &  1.409e-3 & 2.00 & 1.317e-3 & 1.21 \\ \hline
			        1 & 1/32 & 1.891e-6 & 3.00 &  1.418e-6 & 3.00 &  3.517e-4 & 2.00 & 4.188e-4 & 1.65 \\ \hline
			        1 & 1/64 & 2.354e-7 & 3.00 & 1.772e-7 & 3.00 & 8.787e-5 &  2.00 & 1.111e-4 & 1.92  \\ \hline
			      1 & 1/128& 2.938e-8 & 2.99 & 2.229e-8   & 2.99 & 2.196e-5 &  2.00 & 2.835e-5 & 1.97 \\ \hline
			\end{tabular} 
                        \begin{tabular}{|c | c | c | c || c | c || c | c || c | c |}
                \hline
                & & \multicolumn{2}{c ||}{TH $(P_2,P_1)$} & \multicolumn{2}{c || }{SV $(P_2,P_1^{disc})$} 
                & \multicolumn{2}{c ||}{Mini $(P_1^b,P_1)$} & \multicolumn{2}{c | }{CR $(P_1^{nc},P_0)$} 
                \\ \hline 
                $\nu$ & $h$ & $\| (u- u_h)(T) \|_0 $ & Rate & $\| (u- u_h)(T) \|_0 $ & Rate & $\| (u- u_h)(T) \|_0 $ & Rate & $\| (u- u_h)(T) \|_0 $ & Rate  \\ \hline
                $10^{-6}$ &  1/8 & 1.062e-3 & -      & 9.046e-5 & - & 5.448e-3    & -    &      5.294e-3 & - \\ \hline
                $10^{-6}$ & 1/16 & 5.566e-4 & 0.93 & 1.132e-5 &  3.00 &  1.427e-3 & 1.93 & 5.074e-3 & 0.06 \\ \hline
                    $10^{-6}$ & 1/32 & 2.822e-4 & 0.98 &  1.417e-6 & 3.00 &  4.176e-4 & 1.77 & 5.106e-3 & -0.01 \\ \hline
                    $10^{-6}$ & 1/64 & 1.416e-4 & 0.99 &  1.772e-7 & 3.00 & 1.508e-4 &  1.47 & 5.126e-3& 0.01 \\ \hline
                    $10^{-6}$ & 1/128& 7.079e-5 & 1.00 & 2.215e-8   & 3.00 & 6.615e-5 &  1.19 & 5.135e-3 & 0.00 \\ \hline
            \end{tabular}		
	\caption{\label{hconv} $\LL^2$ velocity errors and rates for the Stokes test problem with $\nu=1$ (top) and $\nu=10^{-6}$ (bottom).}
	}
\end{table}

The first test we consider is on $\Omega=(0,1)^2$, with analytical solution
{\small
\begin{align*}
 \uu(x,y,t)  = \langle \cos(y),\ \sin(x) \rangle^T (1+t), \ \ 
 p(x,y,t)  = \sin(x+y).
\end{align*}
}The forcing function $\ff$ is calculated from \eqref{eq:cont:problem} for a given $\nu$, and inhomogeneous Dirichlet boundary conditions are enforced nodally.  To illustrate our theory, we compute on successively refined uniform triangular meshes that are additionally refined with an Alfeld split \cite{FGN18}, and compute with $(\PP_2,P_1)$ Taylor--Hood (TH), $(\PP_2,P_1^{disc})$ Scott--Vogelius (SV), $(\PP_1^{bub},P_1)$ mini, and $(\PP_1^{nc},P_0)$ Crouzeix--Raviart (CR) elements.  To isolate the spatial error, we use BDF3 time stepping with $\Delta t=$1e-3 and end time of $T=0.01$ (using initial conditions taken to be the nodal interpolant of the true solution at $0,\ \Delta t,\ 2\Delta t$).

For each element choice, $\LL^2$ velocity errors and rates are computed for two viscosities, $\nu=1$ and $\nu=10^{-6}$, see table \ref{hconv}.  For $\nu=1$, we observe optimal convergence for all elements
as predicted by the classical theory \cite{HR90}.  However, we observe very different behavior with $\nu=10^{-6}$.  Here, only the pressure-robust Scott--Vogelius elements provide optimal convergence, and all other element choices lose {\it one} (mini element) or
{\it two} convergence orders on preasymptotic meshes,
as is predicted by our novel analysis above.

\subsection{Chorin vortex decay for time dependent Navier-Stokes}

For a second test, we choose the Chorin problem for incompressible NS \cite{Cho68}.  Although our analysis is
for time dependent Stokes, NS is still relevant since the same kind of dominant pressure term exists (however an analysis would be more complex due to the nonlinear term),
since the Chorin problem is a so-called Beltrami flow,
i.e., here the nonlinear term $(\uu \cdot \nabla) \uu=\frac{1}{2} \nabla (|\uu|^2)$ is
a gradient balanced by the pressure gradient.
The domain is taken to be the unit square $\Omega=(0,1)\times(0,1)$, 
and the true NSE solution is taken to be
{\small
\begin{align*}
 \uu(x,y,t)  = \langle -\cos(n\pi x)\sin(n\pi y), \sin(n\pi x)\cos(n\pi y) \rangle^T e^{-2n^{2}\pi^{2}\nu t},\ \ 
 p(x,y,t)  =-\frac{1}{4}(\cos(2n\pi x)+\cos(2n\pi y))e^{-2n^{2}\pi^{2}\nu t},
\end{align*}
} with $n=2$.  This system is an exact solution to the incompressible NS equations with forcing $\ff={\bf 0}$ and $\uu_0=\langle u_1(x,y,0), u_2(x,y,0) \rangle^T$.  We use the same spatial and temporal discretizations as in the first example, and again test with $\nu=1$ and $\nu=10^{-6}$.  Inhomogenous Dirichlet boundary conditions are enforced nodally.  

Results for this test are shown in table \ref{hconv2}, and we observe very similar results to the Stokes test problem
above.  For large $\nu$, all tests show optimal convergence. For $\nu=10^{-6}$, SV
error appears to converge with second order, while CR  error locks, and both mini and TH element solutions converge with just first order.
\begin{table}[h!]
{\footnotesize
			\centering
			
						\begin{tabular}{|c|c | c | c || c | c || c | c || c | c |}
				\hline
				& & \multicolumn{2}{c ||}{TH $(P_2,P_1)$} & \multicolumn{2}{c || }{SV $(P_2,P_1^{disc})$} 
				& \multicolumn{2}{c ||}{Mini $(P_1^b,P_1)$} & \multicolumn{2}{c | }{CR $(P_1^{nc},P_0)$} 
				\\ \hline 
				$\nu$ & $h$ & $\| (u- u_h)(T) \|_0 $ & Rate & $\| (u- u_h)(T) \|_0 $ & Rate & $\| (u- u_h)(T) \|_0 $ & Rate & $\| (u- u_h)(T) \|_0 $ & Rate  \\ \hline
				1 & 1/8 & 1.751e-2 & -      & 4.696e-2 & - & 2.004e-1    & -    &     5.956e-2& - \\ \hline
				1 & 1/16 & 2.203e-2 & 2.99 & 6.475e-3 & 2.85  &  7.683e-2 & 1.38 & 1.607e-2 & 1.89 \\ \hline
			        1 & 1/32 & 2.846e-4 & 2.96 &  8.423e-4& 2.94  &  2.175e-2 & 1.82 & 4.094e-3& 1.97 \\ \hline
			        1 & 1/64 & 3.634e-5 & 2.97& 1.063e-4 & 2.99  & 5.594e-3 &  1.96 & 1.031e-3& 1.99  \\ \hline
			        1 & 1/128& 4.572e-6 & 2.99 & 1.451e-5    & 2.87  & 1.401e-3 &  2.00 & 2.610e-4 & 1.98 \\ \hline
			\end{tabular}			
				\begin{tabular}{|c|c | c | c || c | c || c | c || c | c |}
				\hline
				& & \multicolumn{2}{c ||}{TH $(P_2,P_1)$} & \multicolumn{2}{c || }{SV $(P_2,P_1^{disc})$} 
				& \multicolumn{2}{c ||}{Mini $(P_1^b,P_1)$} & \multicolumn{2}{c | }{CR $(P_1^{nc},P_0)$} 
				\\ \hline 
				$\nu$ & $h$ & $\| (u- u_h)(T) \|_0 $ & Rate & $\| (u- u_h)(T) \|_0 $ & Rate & $\| (u- u_h)(T) \|_0 $ & Rate & $\| (u- u_h)(T) \|_0 $ & Rate  \\ \hline
				$10^{-6}$ & 1/8 & 2.470e-2 & -      & 7.242e-2 & - & 1.310e-1    & -    &      9.357e-2 & - \\ \hline
				$10^{-6}$ & 1/16 & 8.441e-3 & 1.55 & 1.083e-2 &  2.74 &  3.792e-2 & 1.81 & 2.920e-2 & 1.68 \\ \hline
			        $10^{-6}$ & 1/32 & 3.899e-3 & 1.11 &  1.682e-3 & 2.69 &  1.077e-2 & 1.79 & 1.836e-2 & 0.67 \\ \hline
			        $10^{-6}$ & 1/64 & 1.879e-3 & 1.05 & 2.677e-4 & 2.65 & 3.555e-3 &  1.60 & 1.753e-2& 0.07  \\ \hline
			        $10^{-6}$ & 1/128& 8.481e-4 & 1.15 & 5.004e-5   & 2.42  & 1.416e-3 &  1.33 & 1.759e-2 & 0.00 \\ \hline
			\end{tabular}			
	\caption{\label{hconv2} $\LL^2$ velocity errors and rates for the Chorin test with $\nu=1$ (top) and $\nu=10^{-6}$ (bottom).}
	}
\end{table}
\vspace{-0.5cm}
\section{Conclusions}
While it is well known that `optimal' theoretical convergence rates are often not observed when $\nu \ll 1$ except on very fine meshes, little seems known about how error behaves on computable meshes.  We gave herein a new and sharp numerical
analysis for the $\LL^2$ velocity error in the {\em time-dependent} Stokes equations,
emphasizing the role gradient forces for the error evolution.
In particular, two cases arise: if classical (non-pressure-robust) elements are used, suboptimal convergence (by two orders for TH-like element families, or one order by equal-order elements)
and locking will occur, but if pressure-robust elements are used,
optimal $\LL^2$ convergence can be maintained.
Note that no a-priori knowledge of the equilibrium solution
is required,
which is a typical disadvantage of well-balanced
schemes for hyperbolic conservation laws \cite{HAT18};
the $\LL^2$-orthogonality of certain velocity test functions against arbitrary gradient fields suffices \cite{JLMNR17}. \\ \ \\
{\bf Acknowledgements}
The author L.R. acknowledges support from National Science Foundation Grant DMS 1522191.

\bibliographystyle{abbrvurl}
\bibliography{bibliography}

\end{document}